\documentclass[11pt]{article}
\usepackage{authblk}

\usepackage{geometry}
\usepackage{amsmath,amssymb,mathrsfs}
\usepackage{color}
\usepackage{graphicx}
\usepackage{subfigure}
\usepackage{soul}
\usepackage{bbding}
\usepackage{lineno,hyperref}
\usepackage{multirow, booktabs, stmaryrd}
\modulolinenumbers[5]

\def\bp{\mathbf{p}}
\def\bx{\mathbf{x}}

\newcommand{\fr}[2]{\frac{#1}{#2}}

\newcommand{\blue}[1]{\textcolor{black}{#1}}

\title{A Discontinuity Capturing Shallow Neural Network for Elliptic Interface Problems}

\author[1,3]{Wei-Fan Hu}
\author[2,3]{Te-Sheng Lin}
\author[2]{Ming-Chih Lai}

\affil[1]{Department of Mathematics, National Central University, Taoyuan 32001, Taiwan}
\affil[2]{Department of Applied Mathematics, National Yang Ming Chiao Tung University, Hsinchu 30010, Taiwan}
\affil[3]{National Center for Theoretical Sciences, National Taiwan University, Taipei 10617, Taiwan}

\begin{document}

\maketitle
\begin{abstract}
In this paper, a new Discontinuity Capturing Shallow Neural Network (DCSNN) for approximating $d$-dimensional piecewise continuous functions and for solving elliptic interface problems is developed. There are three novel features in the present network; namely, (i) jump discontinuities are accurately captured, (ii) it is completely shallow, comprising only one hidden layer, (iii) it is completely mesh-free for solving partial differential equations. The crucial idea here is that a $d$-dimensional piecewise continuous function can be extended to a continuous function defined in $(d+1)$-dimensional space, where the augmented coordinate variable labels the pieces of each sub-domain. We then construct a shallow neural network to express this new function. Since only one hidden layer is employed, the number of training parameters (weights and biases) scales linearly with the dimension and the neurons used in the hidden layer. For solving elliptic interface problems, the network is trained by minimizing the mean square error loss that consists of the residual of the governing equation, boundary condition, and the interface jump conditions. We perform a series of numerical tests to demonstrate the accuracy of the present network. Our DCSNN model is efficient due to only a moderate number of parameters needed to be trained (a few hundred parameters used throughout all numerical examples), and the results indicate good accuracy. Compared with the results obtained by the traditional grid-based immersed interface method (IIM), which is designed particularly for elliptic interface problems, our network model shows a better accuracy than IIM. We conclude by solving a six-dimensional problem to demonstrate the capability of the present network for high-dimensional applications.
\noindent {\it keywords}: Shallow neural network, deep learning, discontinuity capturing, elliptic interface problem, high-dimensional PDEs
\end{abstract}

\section{Introduction}

For the past decade, deep learning has gained great success in image recognition, natural language processing, computer vision, and other numerous practical applications in our daily life. Until very recently, \blue{it has begun to} draw much attention to solve partial differential equations (PDEs) using deep neural networks (DNNs) \blue{in the scientific computing community}. Part of the theoretical reason can be attributed to the various kinds of expressive power for function approximations using DNN, such as those described in \cite{Cybenko1989, HSW89, LPWHW17, HS18}, and the references therein. The usage of automatic differentiation in machine \blue{learning~\cite{BPRS18}, making} derivatives evaluations through neural networks \blue{much easier, is probably} another reason in practice. However, from the authors' points of view, \blue{the most motivating reason could} be when the PDEs become intractable to traditional numerical methods, such as finite difference, finite element, or spectral method. For instance, even with today's growing computing power and resources, solving high-dimensional PDEs is still notoriously difficult to the above traditional numerical methods due to the ``curse of dimensionality". Some successful methods for solving high-dimensional PDEs using the deep learning approach can be found in \cite{Weinan_2018_2, DGM_2018, ROLNF20, LCLHL22}, just to name a few.

As deep neural networks to solve PDEs become popular, there are mainly two different approaches; namely, the physics-informed neural networks (PINNs) by Raissi et al.~\cite{raissi2019} (also termed as DGM in \cite{DGM_2018}), and the deep Ritz method by E and Yu~\cite{weinan_2018}. Both methods use the neural network approximations for solutions of PDEs and share the mesh-free advantage which the popular traditional numerical methods are unlikely to have. The PINNs solve PDEs by minimizing the mean squared error loss of the equation residual, and the initial and boundary condition errors simultaneously. This kind of neural network approach can be dated back to Dissanayake and Phan-Thien~\cite{DP94}, where the authors used a small number of \blue{training points and a two-hidden-layer network} to solve the linear and nonlinear Poisson equations in a two-dimensional square domain. Shortly after, a superposition of the boundary condition and neural network approximation was developed by Lagaris et al.~\cite{LLF98} to solve ordinary and partial differential equations in regular domains. The extension to solve boundary value problems in irregular domains by the synergy of a single hidden layer network and a radial basis function network (for exact satisfaction of the boundary condition) was also proposed in~\cite{LLP00}. On the other hand, the deep Ritz method begins with solving the variational problem equivalent to the original PDE, so the natural loss function in this framework is simply the energy functional. Recently, a penalty-free Ritz method was proposed in \cite{PFNN_2021} to solve a class of second-order boundary value problems in complex domains.

Despite the success of the deep learning approach in solving PDEs with smooth solutions, not much attention was paid to problems with non-smooth solutions, or problems where the solution is only smooth in a piecewise manner. A typical example is an elliptic interface problem, where the solution and its derivatives have jump discontinuities across the interface. The intrinsic difficulty might be because the activation functions in a deep neural network are in general smooth, such as the sigmoid function, or at least continuous, such as the rectified linear unit (ReLU). Along this line, a deep Ritz type approach to solve the elliptic interface problem with high-contrast discontinuous coefficients was developed in \cite{wang_2020}. The network architecture is like the one proposed in \cite{weinan_2018} where two fully connected hidden layers compose a block, so the network consists of two blocks with an output linear layer. Each block is linked by a residual connection. However, there is a lack of intuitive explanation on why such a network is built. Another deep least squares method~\cite{DLS_2020} was proposed to minimize the first-order system of least-squares functional, which is rewritten from the second-order elliptic problem. An accuracy comparison using different loss functions and activation functions is performed. However, the results are only given for one-dimensional problems. Notice that the solutions in both papers~\cite{wang_2020, DLS_2020} are continuous, but the derivatives have jumps across the interface.

To the best of our knowledge, so far the most successful methodology for solving elliptic problems with discontinuous solution and derivatives across an interface is to use the piecewise deep neural network proposed by He et al.~\cite{he2020meshfree}. They approximate the solution by two neural networks corresponding to two disjoint sub-domains, so that in each sub-domain, the solution remains smooth. These two networks are linked by imposing the solution and its normal derivative jump conditions in the mean squared \blue{error loss. This} way, the resulting error of the training solution can be significantly reduced. However, like most DNNs, as the network becomes deeper, the training process requires more computational effort. Not to mention that one has to train an individual neural network in each sub-domain.

In this paper, we propose a new discontinuity capturing shallow neural network (DCSNN) for solving elliptic interface problems. The novelty of the proposed network is three-fold. Firstly, the network captures the solution and its normal derivative jumps sharply across the interface. Secondly, the network is completely shallow, meaning that only one hidden layer is needed, so it significantly reduces the training cost in contrast to DNN. Lastly, it is totally mesh-free, thus it can be applied to solve problems with complex (or irregular) domains which the traditional accurate immersed interface method (IIM)~\cite{Li2006} is hard to implement, especially in higher dimensions.

The rest of the paper is organized as follows. In Section~2, we show how a $d$-dimensional piecewise continuous function can be approximated by the present shallow network. It demonstrates the main underlying idea of the developed DCSNN by function approximation. Then an optimization algorithm to train the network is introduced, followed by an example to confirm the accuracy. The present network for solving elliptic interface problems is presented in Section~3, followed by a series of numerical accuracy tests and comparisons in Section~4. Some concluding remarks and future work are given in Section~5.

\section{Discontinuity capturing shallow neural network}

In this section, \blue{we start by constructing} a neural network to approximate a $d$-dimensional, \emph{piecewise continuous}, scalar function $\phi(\bx)$ in the domain $\Omega = \Omega^-\cup\Omega^+ \cup \Gamma$ ($\Omega\subseteq\mathbb{R}^d$) defined by
\begin{align}
\phi(\bx) =
\left\{
\begin{array}{ll}
\phi^-(\bx)  & \mbox{if\;\;} \bx \in \Omega^-,\\
\phi^+(\bx) & \mbox{if\;\;} \bx \in \Omega^+,
\end{array}\right.
\end{align}
where $\bx = (x_1,x_2,\cdots,x_d)$, $\phi^-$ and $\phi^+$ are both smooth functions in their corresponding sub-domains. The interface $\Gamma$ is defined as the boundary between the sub-domains $\Omega^-$ and $\Omega^+$, where the function $\phi$ has a jump discontinuity at $\Gamma$. (\blue{Note that the function} value of $\phi$ on the interface $\Gamma$ can be simply defined by $\phi^-$ for convenience purpose.)  In order to approximate the above piecewise continuous function in $d$ dimensions, we first make a continuous extension of the function in $(d+1)$-dimensional space, and then construct a shallow (one hidden layer) neural network to approximate such a function. Since we aim to approximate the piecewise continuous function accurately and capture the function discontinuity sharply using only a one-hidden-layer neural network, we refer to the present network approximation as the Discontinuity Capturing Shallow Neural Network (DCSNN).

\subsection{Continuous function extension}

The key idea \blue{here is that a $d$-dimensional} discontinuous function can be alternatively represented by a $(d+1)$-dimensional smooth function. To this end, we first introduce an augmented one-dimensional variable $z$ to categorize the variable $\bx$ into two types; namely, $\bx$ in $\Omega^-$ or $\Omega^+$. More precisely, we set $z = -1$ if $\bx\in\Omega^-$ and $z = 1$ if $\bx\in\Omega^+$ (this task can be easily done with the help of a level set function for which the zero level set denotes $\Gamma$). We then define the $(d+1)$-dimensional function using the augmentation variable $(\bx,z)$ as
\begin{align}
\phi_{aug}(\bx,z) =
\left\{
\begin{array}{ll}
\phi^-(\bx)  & \mbox{if\;\;} z = -1,\\
\phi^+(\bx) & \mbox{if\;\;} z = 1,
\end{array}\right.
\end{align}
where $\bx\in \Omega$ and $z\in \mathbb{R}$.
That is, both $\phi^-$ and $\phi^+$ are regarded as smooth extensions over the entire domain $\Omega$, and the augmented function $\phi_{aug}(\bx,z)$ is assumed to be continuous on the domain $\Omega\times\mathbb{R}$. In such a way,
 the function $\phi$ can be rewritten in terms of the augmented function as
\begin{align}\label{Eq:augphi}
\phi(\bx) =
\left\{
\begin{array}{ll}
\phi_{aug}(\bx, -1) & \mbox{if\;\;} \bx \in \Omega^-,\\
\phi_{aug}(\bx, 1)  & \mbox{if\;\;} \bx \in \Omega^+.\\
\end{array}\right.
\end{align}
Therefore, the piecewise continuous function $\phi$ now can be regarded as a continuous function defined on one-dimensional higher space and restricted to its sub-domains.

\begin{figure}[h]
\centering
\includegraphics[scale=0.33]{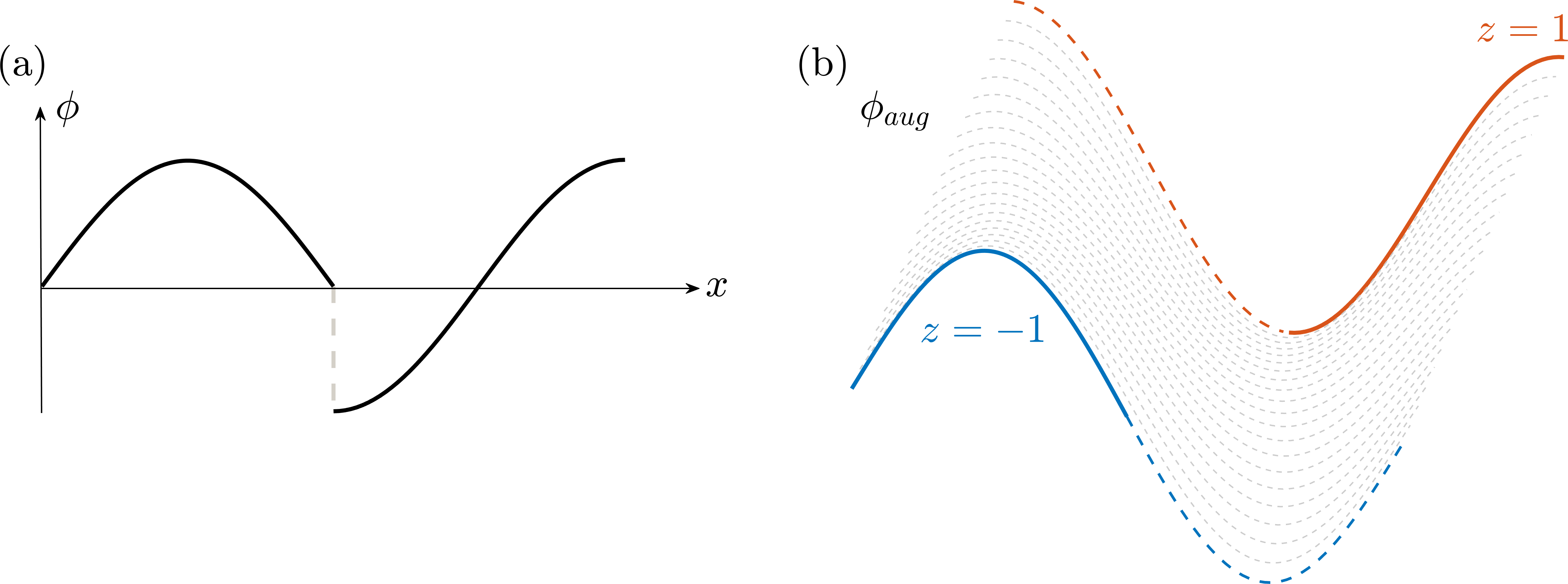}
\caption{(a) The one-dimensional piecewise continuous function $\phi(x)$.
(b) The two-dimensional augmented function $\phi_{aug}(x,z)$ at different $z$ values.
}
\label{Fig:phi_aug}
\end{figure}

We illustrate this idea by considering a one-dimensional example where $\phi^-(x) = \sin(2\pi x)$ if $x\in[0,\frac{1}{2})$ and $\phi^+(x) = \cos(2\pi x)$ if $x\in(\frac{1}{2},1]$. The function has a jump discontinuity at $x=\frac{1}{2}$, see Fig.~\ref{Fig:phi_aug}(a). The two-dimensional augmented function can be constructed as $\phi_{aug}(x,z) = \fr{1-z}{2}\phi^-(x) + \fr{1+z}{2}\phi^+(x)$, see Fig.~\ref{Fig:phi_aug}(b) for the function profile. In particular, its restriction at $z=-1$ corresponding to $\phi^-$ is shown as the blue solid curve, and the one at $z = 1$ corresponding to $\phi^+$ is shown as the red solid curve. It is important to mention that this augmented approach can be straightforwardly applied to piecewise continuous functions with arbitrary many pieces by simply labelling various $z$ values. Besides, the augmented function is not unique, i.e., there exists infinitely many such functions that have their restrictions to be $\phi$. The remaining issues are how to construct the neural network to approximate the augmented function and how to train the network efficiently, which we shall describe in the next subsection.

\subsection{Shallow neural network structure}\label{subsec:SINN}

\begin{figure}[h]
\centering
\includegraphics[scale=0.35]{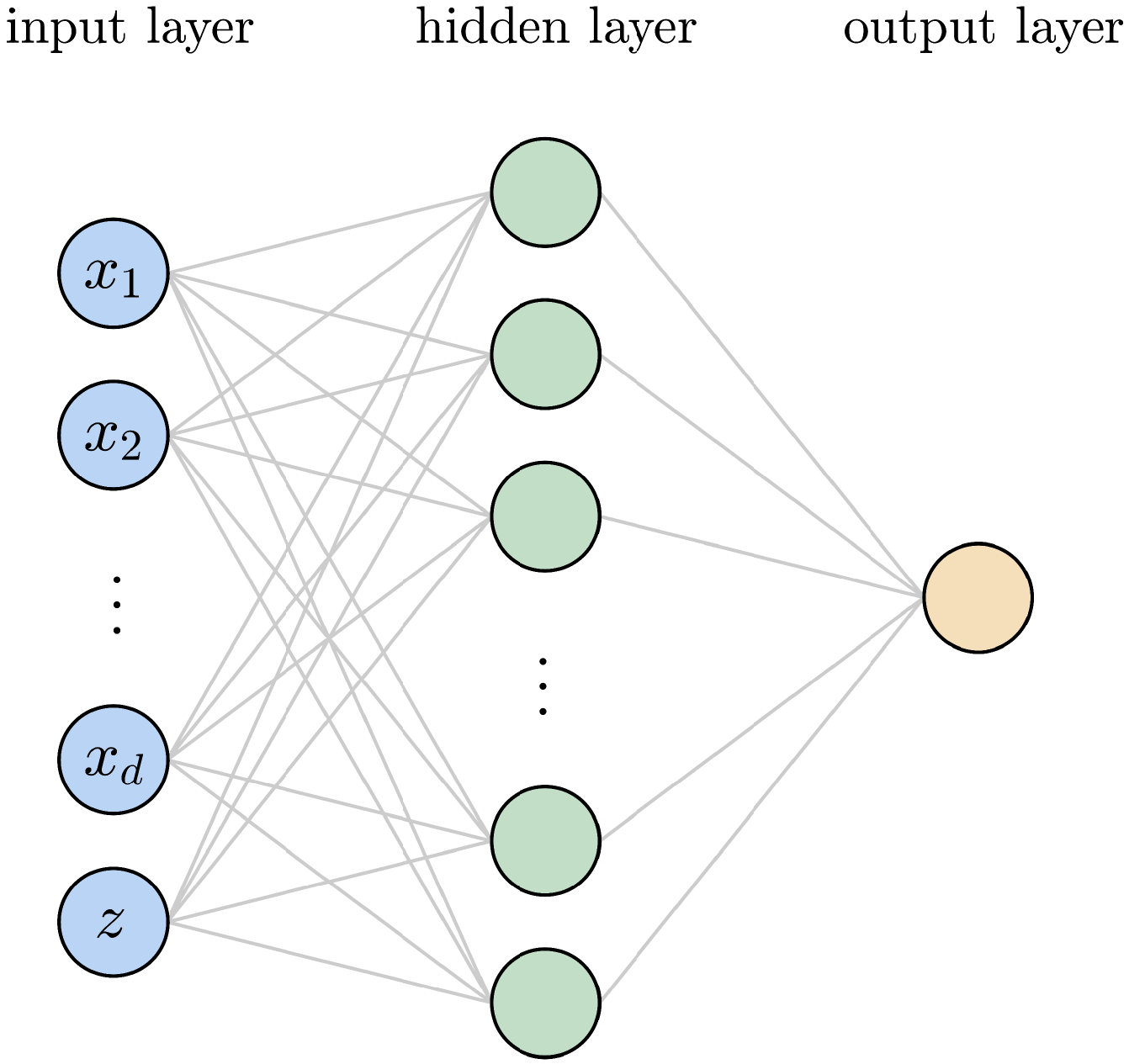}
\caption{The structure of DCSNN. $N$ neurons are used in the hidden layer.}
\label{Fig:SINN}
\end{figure}

We propose a shallow neural network to approximate the $d$-dimensional piecewise continuous function. As shown in the previous subsection, all we have to do is to construct the augmented function $\phi_{aug}$ that is in fact continuous in $d+1$ dimensions. Based on the universal approximation theory~\cite{Cybenko1989}, we hereby design a \emph{shallow}, \emph{feedforward}, \emph{fully connected} neural network architecture, in which only one single hidden layer is employed to approximate $\phi_{aug}$. The structure of this shallow neural network is shown in Fig.~\ref{Fig:SINN}. Let $N$ be the number of neurons used in the hidden layer, the augmented function (or output layer) under this network structure can be explicitly expressed by
\begin{align}\label{Eq:SINN}
\phi_{aug}(\bx,z) = W^{[2]}\sigma(W^{[1]}(\bx,z)^T+b^{[1]}) + b^{[2]},
\end{align}
where $W^{[1]}\in\mathbb{R}^{N\times(d+1)}$ and $W^{[2]}\in\mathbb{R}^{1 \times N}$ \blue{are} the weights, $b^{[1]}\in\mathbb{R}^{N}$ and $b^{[2]}\in\mathbb{R}$ \blue{are} the biases, and $\sigma$ is the activation function. One can easily see that the augmented function is simply a finite linear combination of activation functions.

By collecting all the training parameters (including all the weights and biases) in a vector $\bp$, the total number of parameters in the network (i.e., dimension of $\bp$) is counted by $N_p = (d+3)N + 1$. Given $M$ training data with feature inputs $\{(\bx^i,z^i)\}^M_{i=1}$ and target outputs $\{\phi(\bx^i)\}^M_{i=1}$, where $z^i$ is determined by identifying the category of $\bx^i$, these parameters in $\phi_{aug}$ are learned via minimizing the following mean squared error loss
\begin{align}\label{Eq:loss}
\mbox{Loss}(\bp) = \fr{1}{M}\sum_{i = 1}^M\left( \phi(\bx^i) - \phi_{aug}(\bx^i,z^i;\bp) \right)^2.
\end{align}

Now, we are in a position to solve the above optimization problem. Since the present shallow network generally results in a moderate number of training parameters, the above optimization problem can be efficiently solved by Levenberg-Marquardt (LM) method~\cite{Marquardt63}, a full-batch optimization algorithm. That is, we find $\bp$ via the iterative procedure
\begin{align}\label{Eq:LM}
\bp^{(k+1)} = \bp^{(k)} + \left( J^TJ + \mu I \right)^{-1}\left[ J^T\left(\Phi-\Phi_{aug}(\bp^{(k)})\right) \right],
\end{align}
where the $i$th component of the vectors $\Phi$ and $\Phi_{aug}$ correspond to $\phi(\bx^i)$ and $\phi_{aug}(\bx^i,z^i;\bp)$ in the loss function (\ref{Eq:loss}), and $J\in\mathbb{R}^{M\times N_{p}}$ (typically $M>N_p$) is the Jacobian matrix given by $J = \partial\Phi_{aug}/\partial\bp$. The computation \blue{of the Jacobian matrix} can be done using automatic differentiation \cite{BPRS18, GW08}, or backpropagation algorithm. Notice that $\mu$ is the damping parameter. The \blue{general strategy for tuning the} damping parameter is as follows. Initially, $\mu$ is set to be large \blue{so that the first few updates }are small steps along the gradient direction (note that the term $J^T\left(\Phi-\Phi_{aug}(\bp)\right)$ denotes the negative gradient direction of the loss function). We adjust the parameter $\mu$ lower as the loss decreases, in such a way, the Levenberg-Marquardt method approaches the Gauss-Newton method, and the solution to the optimization problem typically speeds up to the local minimum. The major cost in each iteration step (\ref{Eq:LM}) comes from the computation of matrix-vector multiplication. This can be done efficiently by using reduced singular value decomposition to the Jacobian matrix $J = U\Sigma V^T$ so that the updating direction (the second term on the right-hand side of Eq.~(\ref{Eq:LM})) is computed by $V\mbox{diag}(\fr{\sigma_J}{\sigma_J^2+\mu})U^T(\Phi-\Phi_{aug})$, where $\sigma_J$ denotes the singular values of $J$.

\subsection{An example of function approximation}

For illustration purpose, here we show the capability of the present network by considering a one-dimensional piecewise continuous function
\begin{align}
\phi(x) =
\left\{
\begin{array}{ll}
\sin(2\pi x)  & \mbox{if\;\;} x \in [0,\fr{1}{2}),\\
\cos(2\pi x) & \mbox{if\;\;} x \in (\fr{1}{2},1].
\end{array}\right.
\end{align}
The network is trained by using only $5$ neurons in the hidden layer, thus there are totally $21$ parameters in weights and biases needed to be learned.
We use $100$ randomly sampled training points in the interval $[0,1]$ (including two boundary points) with the sigmoid activation function and stop the iteration when $\mbox{Loss}(\bp) < 10^{-12}$.

\begin{figure}[h]
\centering
\includegraphics[scale=0.4]{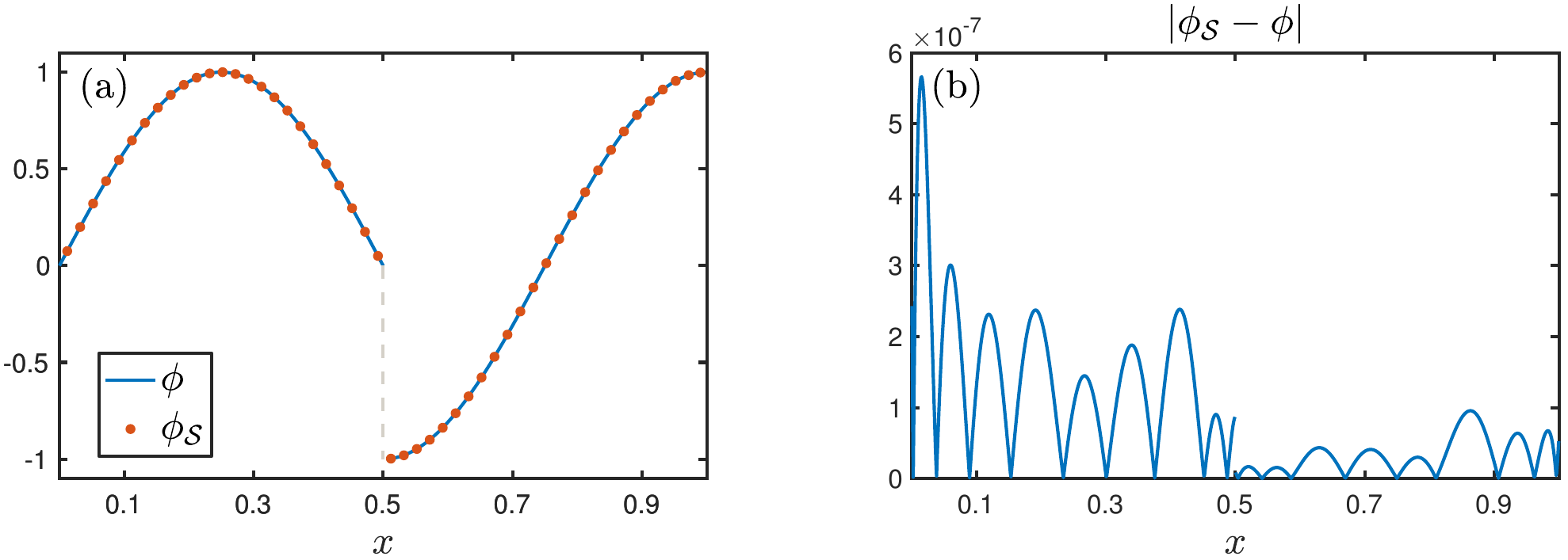}
\caption{(a) The piecewise continuous function $\phi$ (denoted by the solid line) and the DCSNN function $\phi_\mathcal{S}$ (denoted by $50$ dotted points among $1000$ test points). (b) The plot of absolute errors $|\phi_\mathcal{S}-\phi|$ at $1000$ test points. The maximum absolute error is $\|\phi_\mathcal{S}-\phi\|_\infty = 5.661\times10^{-7}$.}
\label{Fig:1d_SINN}
\end{figure}

The results are shown in Fig.~\ref{Fig:1d_SINN}. In panel (a) we compare the exact solution $\phi$ with the DCSNN solution, denoted by $\phi_\mathcal{S}$. Both of them are evaluated at $1000$ test points, which are equally distributed in $[0,1]$. As can be seen, the two functions agree well with each other. It also clearly shows that a jump discontinuity at $x=0.5$ is captured sharply. We show in panel (b) the pointwise absolute error $|\phi_\mathcal{S}-\phi|$. It is interesting to see that the maximum error occurs close to the domain boundary rather than at the point of discontinuity. The prediction of the present network is quite accurate where the maximum norm of the error is of magnitude $O(10^{-7})$. It is important to mention that, despite only one-dimensional case presented here, discontinuous functions in higher dimensions can also be approximated accurately and efficiently using the present network structure (see the last example in Section 4).

\section{Elliptic interface problems}

Now we apply the DCSNN to solve $d$-dimensional elliptic interface problems with piecewise constant coefficients. In such a problem, a domain $\Omega \subset \mathbb{R}^d$ is divided by an embedded interface $\Gamma \subset \mathbb{R}^{d-1}$ into two regions, namely, inside ($\Omega^-$) and outside ($\Omega^+$) of the interface, so that $\Omega = \Omega^- \cup \Omega^+ \cup \Gamma$.
The elliptic interface equation with general inhomogeneous jump conditions reads
\begin{align}\label{Eq:elliptic}
\nabla\cdot(\beta\nabla\phi) = f \quad \mbox{in }\Omega \setminus \Gamma, \quad [\phi]=v, \quad [\beta\partial_n\phi]=w, \quad \mbox{on }\Gamma,
\end{align}
where $\beta$ is a piecewise positive constant function defined by $\beta = \beta^-$ in $\Omega^-$ and $\beta = \beta^+$ in $\Omega^+$. (Note that the present method has no difficulty for handling the variable coefficient case, i.e., $\beta$ is spatial dependent, as one can see from the implementation later.) Here the bracket $[\cdot]$ denotes the jump discontinuity of the quantity approaching from $\Omega^+$ minus the one from $\Omega^-$, $\mathbf{n}$ stands for the outward normal vector defined on the interface $\Gamma$ and $\partial_n\phi$ is the shorthand notation for $\nabla\phi\cdot\mathbf{n}$. Certainly, the above differential equation should be accompanied with suitable boundary conditions (say Dirichlet or Neumann) along the domain boundary $\partial\Omega$. Throughout this paper we will focus on the case of the Dirichlet boundary condition $\phi|_{\partial\Omega} = g$, while other types of boundary conditions will not change the main ingredients presented here.

We first rewrite Eq.~(\ref{Eq:elliptic}) in the form of the Poisson equation:
\begin{align}\label{Eq:Poisson}
\Delta\phi = \tilde{f} =
\left\{
\begin{array}{ll}
f^-/\beta^-   & \mbox{if\;\;} \bx\in \Omega^-\\
f^+/\beta^+ & \mbox{if\;\;} \bx \in \Omega^+\\
\end{array}\right.,  \quad  [\phi]=v, \quad [\beta\partial_n\phi]=w, \quad \mbox{on }\Gamma,
\end{align}
and introduce a $(d+1)$-dimensional augmented function $\phi_{aug}(\bx, z)$ that satisfies
\begin{align}\label{Eq:Poisson_aug}
\Delta_{\bx}\phi_{aug}(\bx, z)  =
\left\{
\begin{array}{ll}
f^-/\beta^- & \mbox{if\;\;} \bx\in\Omega^-,\, z = -1\\
f^+/\beta^+ & \mbox{if\;\;} \bx\in\Omega^+,\, z = 1\\
\end{array}\right.,  \quad  [\phi_{aug}] =v, \quad [\beta\partial_n\phi_{aug}]=w, \quad \mbox{on }\Gamma,
\end{align}
where the jump condition is evaluated by taking the difference between values at $z=1$ and $z=-1$, i.e., $[\phi_{aug}]= \phi_{aug}(\bx, 1) - \phi_{aug}(\bx, -1)$ for $\bx\in\Gamma$, and the same manner applies for $[\beta\partial_n\phi_{aug}]$ (note that $\partial_n\phi_{aug} = \nabla_{\bx}\phi_{aug}\cdot\mathbf{n}$). Once $\phi_{aug}$ is found, the function $\phi$ is recovered using Eq.~(\ref{Eq:augphi}).

We then need to solve the augmented elliptic interface problem (\ref{Eq:Poisson_aug}) to obtain $\phi_{aug}$ using the present DCSNN. Given training points in $\Omega$, on the domain boundary $\partial\Omega$ and along the embedded interface $\Gamma$, denoted by $\{(\mathbf{x}^i, z^i)\}_{i=1}^M$, $\{\mathbf{x}^j_{\partial\Omega}\}_{j=1}^{M_b}$ and $\{\mathbf{x}^k_{\Gamma}\}_{k=1}^{M_\Gamma}$, respectively, we hereby solve the equation (\ref{Eq:Poisson_aug}) by minimizing the mean squared error loss in the framework of the physics-informed learning technique~\cite{raissi2019}
\begin{align}\label{Eq:loss_PINN}
\begin{split}
\text{Loss}(\bp) & =  \frac{1}{M}\sum_{i = 1}^M \left( \Delta_\bx\phi_{aug}(\mathbf{x}^i,z^i) - \tilde{f}(\mathbf{x}^i) \right)^2
+ \frac{\alpha_b}{M_b}\sum_{j = 1}^{M_b} \left( \phi_{aug}(\mathbf{x}^j_{\partial\Omega},1) - g(\mathbf{x}^j_{\partial\Omega}) \right)^2\\
& + \frac{\alpha_\Gamma}{M_{\Gamma}}\sum_{k = 1}^{M_\Gamma} \left( [\phi_{aug}(\mathbf{x}^k_\Gamma)] - v(\mathbf{x}^k_\Gamma) \right)^2 + \frac{\alpha_\Gamma}{M_{\Gamma}}\sum_{k = 1}^{M_\Gamma} \left( [\beta\partial_n\phi_{aug}(\mathbf{x}^k_\Gamma)] - w(\mathbf{x}^k_\Gamma) \right)^2,
\end{split}
\end{align}
where $\alpha_b$ and $\alpha_\Gamma$ are positive constants that can be adjusted to strengthen or weaken constraints at the domain boundary and on the interface, respectively. Notice that, in the \blue{above equation, we have} collected all training parameters (weights and biases) in the vector $\bp$. Since the least squares loss is adopted, again, it could be quite efficient to train the network using the Levenberg-Marquardt method as introduced in \ref{subsec:SINN}. We mention that the spatial derivatives of the target function $\phi_{aug}$ in the loss function, Eq.~(\ref{Eq:loss_PINN}), can be computed by automatic differentiation easily. Although the present network is similar in spirit to PINNs~\cite{raissi2019}, here we only use one hidden layer with sufficiently small number of neurons so it reduces the computational complexity and learning workload significantly without sacrificing the accuracy. Moreover, in solving Eq.~(\ref{Eq:elliptic}), as shown in the next section, the present DCSNN not only achieves better accuracy than the traditional finite difference method, such as the immersed interface method~\cite{Li2006, Hu2015, Hsu2019}, but also outperforms other piecewise DNN~\cite{he2020meshfree} in terms of accuracy and network complexity.

\section{Numerical results}
In this section, we use the developed network to perform numerical tests for the elliptic interface problems in two, three, and even six dimensions. Throughout all the examples (except Example~2, in which we compare our results with existing literature), we set $\beta^- = 1$ and $\beta^+ = 10^{-3}$ (so $\beta^-/\beta^+  = 10^3$, the problem with high contrast coefficients) and the penalty parameters $\alpha_b=\alpha_\Gamma = 1$. The exact solution with a $d$-dimensional variable $\bx = (x_1,x_2,\cdots,x_d)$ is given by
\begin{align}
\phi(\bx) =
\left\{
\begin{array}{ll}
\prod\limits_{i=1}^d \exp(x_i) & \mbox{if\;\;} \bx \in \Omega^-,\\
\prod\limits_{i=1}^d \sin(x_i)  & \mbox{if\;\;} \bx \in \Omega^+,\\
\end{array}\right.
\end{align}
so the corresponding right-hand side function $\tilde{f}$ and the jump information $v$ and $w$ in Eq.~(\ref{Eq:Poisson}) can be easily obtained through the above solution.

In the following examples, we choose sigmoid as the activation function and we finish the training process when the stopping condition $\mbox{Loss}(\bp) < \varepsilon$ is met ($\varepsilon$ is set to be at least smaller than $10^{-9}$). We measure the accuracy of the solution using the test error instead of the training error. Precisely, we randomly choose test points $\{\bx^i\}^{N_{test}}_{i=1}\subset\Omega$ to compute the $L^\infty$ and $L^2$ error as
\begin{align*}
\|\phi_\mathcal{S}-\phi\|_\infty = \max_{1\leq i \leq N_{test}}|\phi_\mathcal{S}(\bx^i)-\phi(\bx^i)|, \quad
\|\phi_\mathcal{S}-\phi\|_2 = \sqrt{\frac{1}{N_{test}}\sum_{i=1}^{N_{test}}(\phi_\mathcal{S}(\bx^i)-\phi(\bx^i))^2},
\end{align*}
where $\phi_\mathcal{S}$ is the solution obtained by the present DCSNN. Throughout this paper, we set the number of test points $N_{test}=100M$, with $M$ training points. We also repeat the numerical experiments $10$ times so the test error reported here is an averaged one. The source codes accompanying this manuscript are available on GitHub~\cite{github}.

\paragraph{\textbf{Example 1}}

We first consider a two-dimensional problem with a regular square domain $\Omega = [-1, 1]\times [-1, 1]$ and an ellipse-shaped interface $\Gamma:\left(\displaystyle\frac{x_1}{0.2}\right)^2+\left(\displaystyle\frac{x_2}{0.5}\right)^2=1$. We use $128$ training points ($M=64$ interior points in the computational domain, $M_b=32$ points on the boundary, and $M_\Gamma=32$ points on the interface).

\begin{figure}[h]
\centering
\includegraphics[scale=0.38]{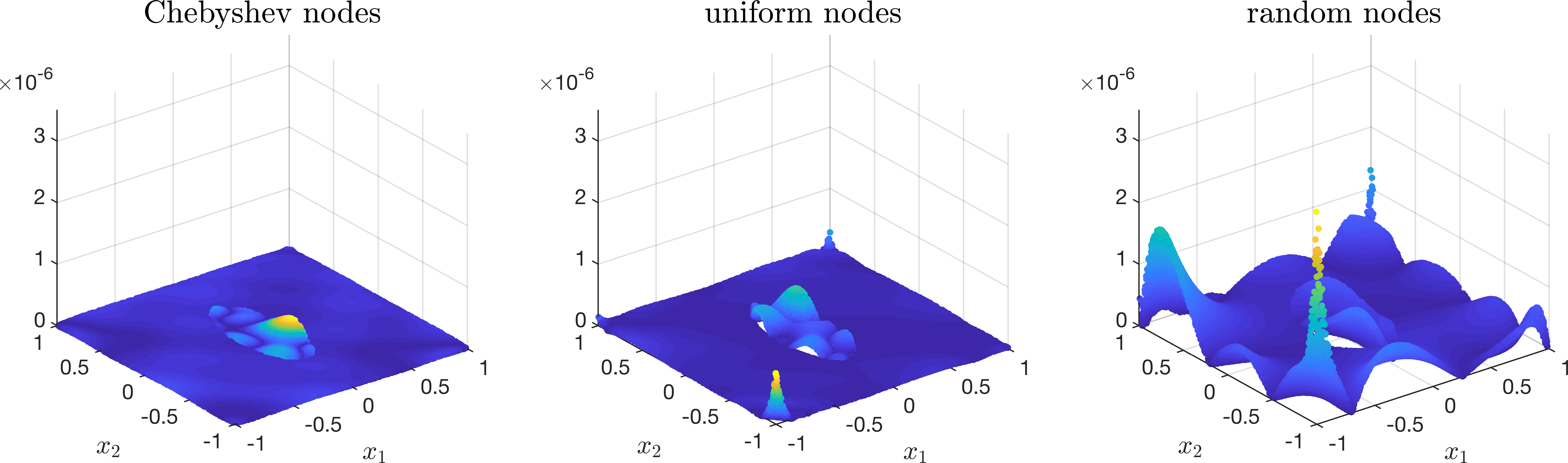}
\caption{The absolute error $|\phi_\mathcal{S}-\phi|$ using $N = 20$ with different training point distributions in Example~1.}
\label{Fig:EX1_error}
\end{figure}

It is interesting to see that the deployment of training data can indeed affect the accuracy of DCSNN solutions. We consider three different distributions of training points; namely, Chebyshev points of the first kind~\cite{Trefethen2000}, uniformly distributed points, and randomly distributed points. The training points on $\Omega$ and $\partial\Omega$ are chosen based on those three distributions, while the points on the interface $\Gamma$ are always chosen randomly. The results are shown in Table~\ref{Table:exp1_node}, where two sets of neurons $N=10, 20$ are used for comparison. We find the errors obtained by Chebyshev and uniform nodes are of the same order of magnitude, whereas the ones obtained by random distribution are less accurate. The corresponding pointwise absolute errors are shown in Fig.~\ref{Fig:EX1_error}. As one can see, the major error of Chebyshev nodes comes from the vicinity of the interface, while the other two distributions have their maximum error occurring close to the corners of the domain. This finding is surprising in the sense that, even though the activation function is fairly different from polynomials, the approximation still achieves better accuracy in Chebyshev nodes.

\begin{table}[h]
\centering
\begin{tabular}{@{}c@{\hspace{10pt}}c@{\hspace{10pt}}c@{\hspace{5pt}}c@{\hspace{10pt}}c@{\hspace{10pt}}c@{\hspace{5pt}}c@{\hspace{10pt}}c@{\hspace{10pt}}c@{\hspace{5pt}}}
\hline
\multirow{2}{*}[-3pt]{$(N,\, N_p)$} & \multicolumn{2}{c}{Chebyshev nodes} & & \multicolumn{2}{c}{uniform nodes} & & \multicolumn{2}{c}{random nodes}\\
\cmidrule(){2-3}\cmidrule(){5-6}\cmidrule(){8-9}
    & $\|\phi_\mathcal{S}-\phi\|_\infty$ & $\|\phi_\mathcal{S}-\phi\|_2$ & & $\|\phi_\mathcal{S}-\phi\|_\infty$ & $\|\phi_\mathcal{S}-\phi\|_2$ & & $\|\phi_\mathcal{S}-\phi\|_\infty$ & $\|\phi_\mathcal{S}-\phi\|_2$\\
\midrule
(10, 51)   & 1.032E$-$05 & 1.168E$-$06 & & 5.965E$-$05 & 8.098E$-$06 & & 3.395E$-$04 & 4.852E$-$05 \\
(20, 101) & 4.386E$-$06 & 7.285E$-$07 & & 8.798E$-$06 & 1.247E$-$06 & & 4.004E$-$05 & 1.574E$-$06 \\
\hline
\end{tabular}
\caption{The errors with different training point distributions in Example~1.}
\label{Table:exp1_node}
\end{table}

\begin{figure}[h]
\centering
\includegraphics[scale=0.38]{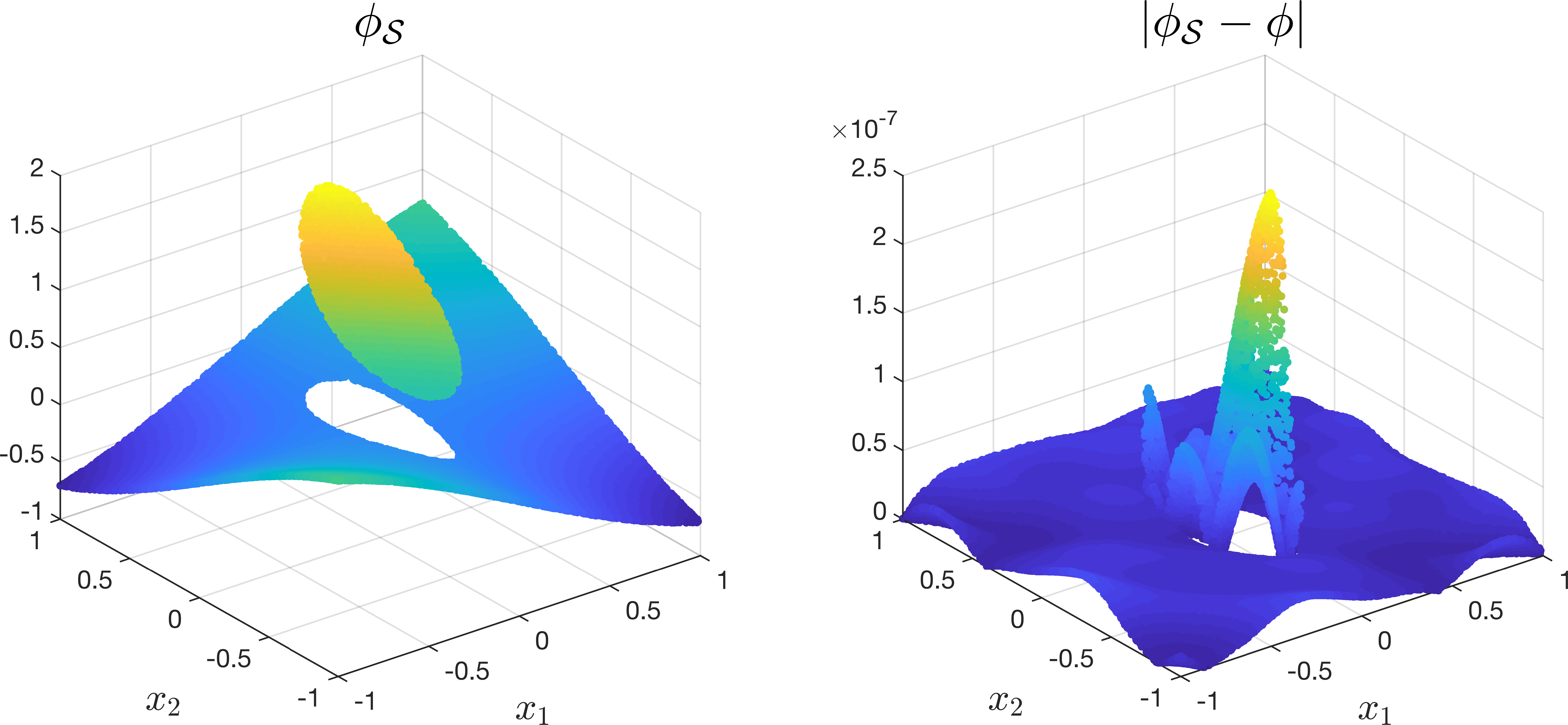}
\caption{Left: The DCSNN solution profile. Right: Absolute error between the DCSNN solution and the exact solution. $N=20$ with Chebyshev training points in Example~1.}
\label{Fig:EX1}
\end{figure}

We show in Fig.~\ref{Fig:EX1} the solution profile of DCSNN; it sharply captures the discontinuity. Meanwhile, it is interesting to see that the absolute error between the exact and the DCSNN solutions also shows a discontinuity across the interface because of the intrinsic construction of the present network.

It is worth mentioning that the number of training points in each sub-domain, $\Omega^-$ and $\Omega^+$, depends roughly on the area ratio between them if the solutions in both sub-domains are equally smooth. \blue{In this example, $|\Omega^-|/|\Omega^+|\sim 0.085$ leads to less than $10$ points in $\Omega^-$.} Even with such a small amount of data, we still get accurate predictions in the interior region. \blue{This might be due to having enough training points on the interface, which brings in an adequate amount of information. Also, note that the number of points needed} in each sub-domain is problem-dependent and often requires some prior information about the solution smoothness in each sub-domain.

We then compare the solution from the Chebyshev nodes (total $128$ training points) in Table~\ref{Table:exp1_node} with the numerical solution obtained by the immersed interface method (IIM)~\cite{Hu2015} which is second-order accurate for the elliptic interface problem Eq.~(\ref{Eq:elliptic}) on Cartesian grids. In IIM, the total number of degree of freedom (number of unknowns) $N_{deg}$ equals the sum of the number of Cartesian grid points $m^2$ and the augmented \blue{projection feet on the interface, $m_\Gamma$; see \cite{Hu2015}.} Table~\ref{Table:exp1} shows the comparison results where the IIM uses the grid resolutions $m=m_\Gamma=128$ and $m=m_\Gamma=256$, so correspondingly $N_{deg}=16512$ and $N_{deg}=65792$, while the number of parameters used in DCSNN are just $N_p=51$ and $N_p=101$, respectively. One can see how significantly different those numbers of unknowns are. With just a few number of neurons, $N=10$ and $20$, the solutions of the present network achieve better accuracy than the IIM. \blue{Similarly, an increase in the number of neurons leads to better accuracy, as expected.}

\begin{table}[h]
\centering
\begin{tabular}{c|c||c|cc}
\hline
 $N_{deg}$ & $\|\phi_{IIM}-\phi\|_{\infty}$ & $(N,\, N_p)$ & $\|\phi_\mathcal{S}-\phi\|_{\infty}$ & $\|\phi_\mathcal{S}-\phi\|_2$ \\
\hline
16512& 7.719E$-$05  & (10,\,51)   & 1.032E$-$05 & 1.168E$-$06 \\
65792& 8.347E$-$06  & (20,\,101) & 4.386E$-$06 & 7.285E$-$07\\
\hline
\end{tabular}
\caption{$\phi$: Exact solution. $\phi_{IIM}$: Solution obtained by IIM. $N_{deg}=16512$ and $65792$ correspond to $m=m_\Gamma=128$ and $m=m_\Gamma=256$. $\phi_\mathcal{S}$: Solution obtained by DCSNN model.}
\label{Table:exp1}
\end{table}

\begin{figure}[h]
\centering
\includegraphics[scale=0.4]{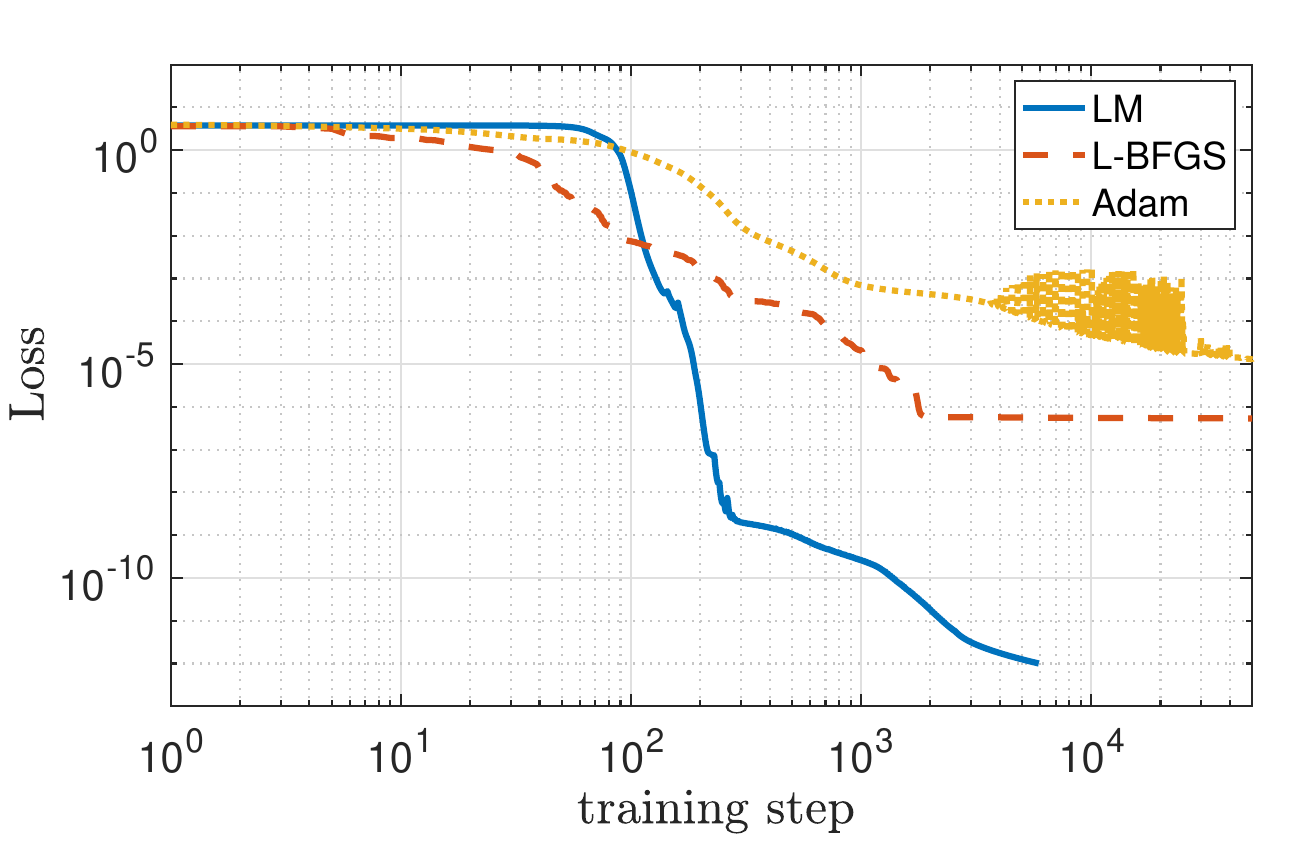}
\caption{The evolution of the loss with different optimizers in Example~1. $N = 20$.}
\label{Fig:EX5_optim}
\end{figure}

One may wonder what will happen if different optimizers are chosen. To see this, we use the same network model, same Chebyshev training points, and solve the same problem, but with Limited-memory Broyden-Fletcher-Goldfarb-Shanno (L-BFGS)~\cite{ZBLN97} and Adaptive Moment Estimation (Adam)~\cite{KB17} optimizers that are commonly used in machine learning practice. We use a learning rate of $0.1$ for L-BFGS, while for Adam we use an initial learning rate $0.01$ and a learning rate scheduler that reduces its value to $0.0001$ during training. We fix the number of iterations to $50000$. The results are shown in Table~\ref{Table:exp1_optimizer}. Compared \blue{with the results from the LM optimizer} in Table~\ref{Table:exp1_node}, the accuracy is about two orders of magnitude lower. We also show the evolution of the loss in Fig.~\ref{Fig:EX5_optim}. The LM loss drops to $10^{-9}$ in about $1000$ steps, while the L-BFGS and Adam losses can not fall below $10^{-7}$ even up to $50000$ training steps.

\begin{table}[h]
\centering
\begin{tabular}{@{}c@{\hspace{10pt}}c@{\hspace{10pt}}c@{\hspace{5pt}}c@{\hspace{10pt}}c@{\hspace{10pt}}c@{\hspace{5pt}}c@{\hspace{10pt}}c@{\hspace{10pt}}c@{\hspace{5pt}}}
\hline
\multirow{2}{*}[-3pt]{$(N,\, N_p)$} &  & \multicolumn{2}{c}{L-BFGS} & & \multicolumn{2}{c}{Adam}\\
\cmidrule(){3-4}\cmidrule(){6-7}
    & & $\|\phi_\mathcal{S}-\phi\|_\infty$ & $\|\phi_\mathcal{S}-\phi\|_2$ & & $\|\phi_\mathcal{S}-\phi\|_\infty$ & $\|\phi_\mathcal{S}-\phi\|_2$\\
\midrule
(10, 51)    & & 1.897E$-$03 & 5.461E$-$04 & & 3.270E$-$03 & 9.077E$-$04 \\
(20, 101) & & 6.086E$-$04 & 1.420E$-$04 & & 1.508E$-$03 & 2.619E$-$04 \\
\hline
\end{tabular}
\caption{Comparison between optimizers in Example~1 with Chebyshev training points.}
\label{Table:exp1_optimizer}
\end{table}

\paragraph{\textbf{Example 2}}

As the second example, we compare the results of DCSNN with the piecewise deep neural network proposed in \cite{he2020meshfree}, in which two (or multiple) individual neural nets are trained to approximate the function in each sub-domain. Here, the square domain $\Omega = [-1,1]^2$ is separated by the interface that is given by the polar curve $r(\theta) = 1/2 + \sin(5\theta)/7$. We choose the exact solution as in \cite{he2020meshfree}
\begin{align*}
\phi(x_1,x_2) =
\left\{
\begin{array}{ll}
\exp(x_1^2+x_2^2) & \mbox{if\;\;} (x_1,x_2)\in\Omega^-,\\
0.1(x_1^2+x_2^2)^2-0.01\log(2\sqrt{x_1^2+x_2^2}) &  \mbox{if\;\;} (x_1,x_2)\in\Omega^+,\\
\end{array}\right.
\end{align*}
and set $\beta^- = 10$ and $\beta^+ = 1$. The present network is trained using $N = 50$ and $100$ neurons with randomly sampled training points $(M,M_b,M_\Gamma) = (400, 80, 80)$. We measure the accuracy of the solutions \blue{using the relative} $L^2$ error and compare the results with piecewise \blue{deep neural nets employing four} and six \blue{hidden layers, as in \cite{he2020meshfree},} for which the same order of training data is used. As one can see from Table~\ref{Table:exp2}, the DCSNN model with only one hidden layer and fewer parameters shows better accuracy.

\begin{table}[h]
\centering
\begin{tabular}{c|c||c|c}
\hline
 $(N,\,N_p)$ & $\|\phi_\mathcal{S}-\phi\|_2/\|\phi\|_2$ & $N_p$ & $\|\phi_\mathcal{DNN}-\phi\|_2/\|\phi\|_2$  \\
\hline
(50,\,251)   &  8.362E$-$04 & 25474 & 4.960E$-$03 \\
(100,\,501) &  2.634E$-$04 & 42114 & 3.727E$-$04 \\
\hline
\end{tabular}
\caption{Relative $L^2$ error for the solutions of DCSNN, $\phi_\mathcal{S}$, and for the solutions of piecewise deep neural network~\cite{he2020meshfree}, $\phi_\mathcal{DNN}$.}
\label{Table:exp2}
\end{table}

\paragraph{\textbf{Example 3}}

\begin{figure}[h]
\centering
\includegraphics[scale=0.38]{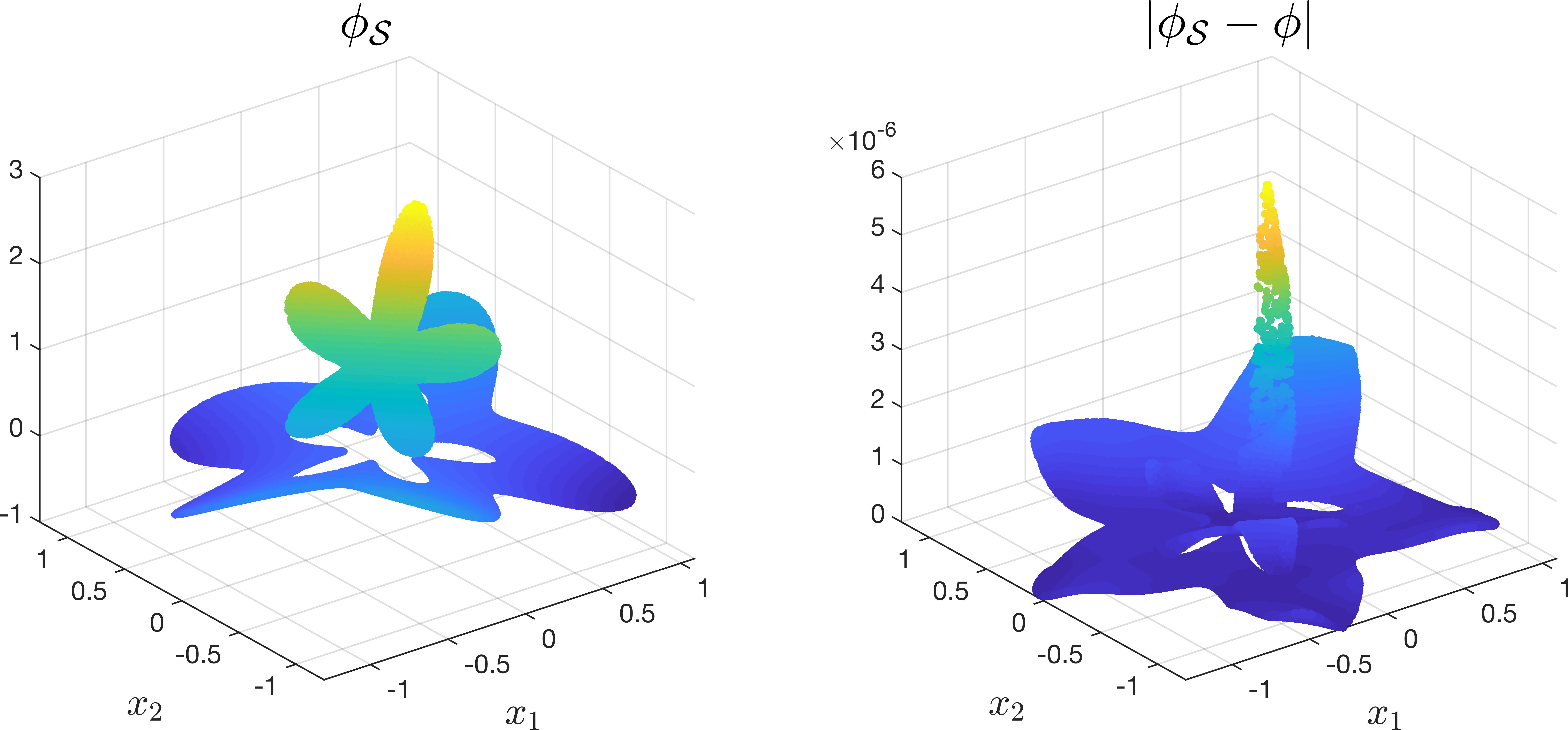}
\caption{Left: The DCSNN solution profile. Right: Absolute error between the DCSNN solution and the exact solution. In this case $\|\phi_\mathcal{S}-\phi\|_\infty = 5.309\times10^{-6}$ and $\|\phi_\mathcal{S}-\phi\|_2=5.856\times10^{-7}$.}
\label{Fig:EX3}
\end{figure}

This example aims to highlight the robustness of DCSNN for handling the elliptic interface problem (\ref{Eq:elliptic}) on irregular domains. Since the proposed method is mesh-free (i.e., not constrained by the locations of training points), the implementation is indeed straightforward. Unlike the previous examples, here the domain is set to be irregular and enclosed by the polar curve $r(\theta) = 1 - 0.3\cos(5\theta)$. The embedded interface $\Gamma$ is also described by the polar curve, $r(\theta) = 0.4 - 0.2\cos(5\theta)$. We use $20$ neurons in the hidden layer (thus $N_p=101$) and randomly \blue{sample} training points with $(M, M_b, M_\Gamma) = (64, 32, 32)$. The results are shown in Fig.~\ref{Fig:EX3}. Without paying extra numerical efforts, the DCSNN model can accurately predict the solution. On the contrary, it can be difficult for traditional finite difference methods, such as IIM, to solve in such a domain. Here, we emphasize that the irregular domains in arbitrary dimensions can be handled properly with no substantial difficulty.

\paragraph{\textbf{Example 4}}

Next, we proceed to the three-dimensional elliptic interface problem. The interface is elliptical, centered at the origin with semi-principal axes of length $0.7$, $0.5$ and $0.3$, embedded in a regular cube $\Omega = [-1,1]^3$. Similar to Example 1, we compare the accuracy among three different training point distributions. The number of training points is $(M,M_b,M_\Gamma) = (216,216,108)$ which is about $4$ times as many as used in the 2D case, while the number of neurons are set by $N=20, 30$ in the hidden layer. The results are shown in Table~\ref{Table:exp4_node}. As expected, the Chebyshev nodes again give the most accurate results compared to the other two.

\begin{table}[h]
\centering
\begin{tabular}{@{}c@{\hspace{10pt}}c@{\hspace{10pt}}c@{\hspace{5pt}}c@{\hspace{10pt}}c@{\hspace{10pt}}c@{\hspace{5pt}}c@{\hspace{10pt}}c@{\hspace{10pt}}c@{\hspace{5pt}}}
\hline
\multirow{2}{*}[-3pt]{$(N, \,N_p)$} & \multicolumn{2}{c}{Chebyshev nodes} & & \multicolumn{2}{c}{uniform nodes} & & \multicolumn{2}{c}{random nodes}\\
\cmidrule(){2-3}\cmidrule(){5-6}\cmidrule(){8-9}
    & $\|\phi_\mathcal{S}-\phi\|_\infty$ & $\|\phi_\mathcal{S}-\phi\|_2$ & & $\|\phi_\mathcal{S}-\phi\|_\infty$ & $\|\phi_\mathcal{S}-\phi\|_2$ & & $\|\phi_\mathcal{S}-\phi\|_\infty$ & $\|\phi_\mathcal{S}-\phi\|_2$\\
\midrule
(20,\,121) & 3.276E$-$05 & 3.832E$-$06 & & 4.214E$-$05 & 2.278E$-$06 & & 1.610E$-$04 & 1.341E$-$05 \\
(30,\,181) & 5.605E$-$06 & 5.809E$-$07 & & 1.733E$-$05 & 1.244E$-$06 & & 7.646E$-$05 & 1.693E$-$06 \\
\hline
\end{tabular}
\caption{The errors with different training point distributions in Example~4.}
\label{Table:exp4_node}
\end{table}

We then compare the accuracy of the DCSNN solution with the 3D immersed interface solver proposed in \cite{Hsu2019}. Note that, in 3D IIM, the total degree of freedom, $N_{deg}$, is now the sum of the number of Cartesian grid points $m^3$ and the augmented projection feet on the interface, $m_\Gamma$. As shown in Table~\ref{Table:exp4}, the results obtained by the present model are more accurate than the ones by IIM, whereas the number of parameters of the network is significantly less than the one used in IIM.

\begin{table}[h]
\centering
\begin{tabular}{c|c||c|cc}
\hline
 $N_{deg}$ & $\|\phi_{IIM}-\phi\|_{\infty}$ & $(N,\, N_p)$ & $\|\phi_\mathcal{S}-\phi\|_{\infty}$ & $\|\phi_\mathcal{S}-\phi\|_2$ \\
 \hline
2201776   & 9.380E$-$04  &  (20,\,121) & 3.276E$-$05 & 3.832E$-$06 \\
17194216 & 2.890E$-$04  & (30,\,181)  & 5.605E$-$06 & 5.809E$-$07 \\
 \hline
\end{tabular}
\caption{$\phi$: Exact solution. $\phi_{IIM}$: Solution obtained by IIM.
$N_{deg}=2201776$ and $17194216$ correspond to $(m,m_\Gamma)=(128,104624)$ and $(m,m_\Gamma)=(256,417000)$.
$\phi_\mathcal{S}$: Solution obtained from DCSNN model.
}
\label{Table:exp4}
\end{table}

It is important to point out that, in DCSNN, the number of neurons in the input layer is the problem dimension plus one. So with a fixed number of neurons in the hidden layer $N$, the total number of parameters $N_p$ increases linearly with the dimension $d$ (recall that $N_p=(d+3)N+1$). In contrast, when the dimensionality of the problem increases, the computational complexity of traditional numerical methods increases significantly.

\paragraph{\textbf{Example 5}}

\begin{figure}[h]
\centering
\includegraphics[scale=0.36]{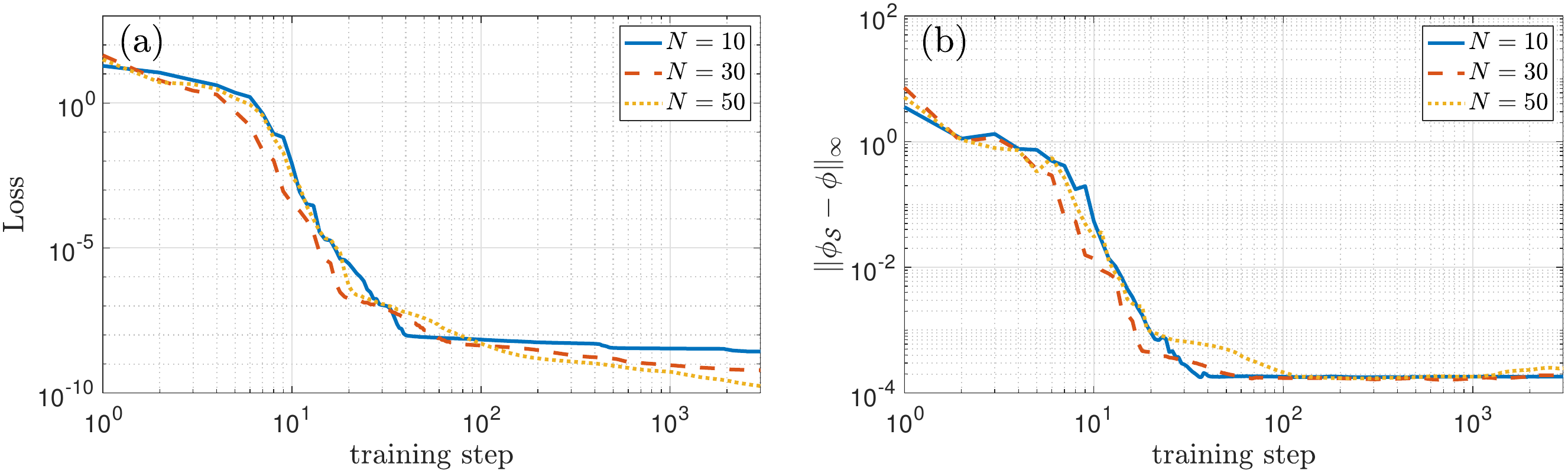}
\caption{The evolutions of (a) the loss and (b) the $L^\infty$ error with $N = 10, 30$ and $50$.}
\label{Fig:EX5}
\end{figure}

As a final example, we show the ability of DCSNN to solve the high-dimensional elliptic interface problem (\ref{Eq:elliptic}) by taking the dimension $d=6$. For the problem setup, we consider a $6$-sphere of radius $0.6$ as the domain $\Omega$ enclosing another smaller $6$-sphere of radius $0.5$ as the interior region $\Omega^-$. We randomly sample the training points with $(M,M_b,M_\Gamma) = (100, 141, 141)$ while the number of neurons are set by $N=10$, $30$ and $50$ (corresponding to total number of parameters $N_p = 91, 271$ and $451$). The evolutions of training loss and $L^\infty$ error against training steps are shown in Fig.~\ref{Fig:EX5}. We find that increasing the number of neurons $N$ leads to a smaller mean squared error loss (see panel (a)). \blue{However, for all cases, the error plateaus in the early training state (within $100$ iterations) and the descent of the loss metric becomes sluggish for the following training steps; see panel (b).} Although not shown here, our experiments show that adding more training points does not improve the accuracy of the predicted solution. We attribute this slow convergence to high dimensionality. \blue{ A similar observation was found} in deep Ritz method for solving high dimensional problems~\cite{weinan_2018}. Despite that, our present solution still gives the error $\|\phi_\mathcal{S}-\phi\|_\infty = 1.689\times10^{-4}$ for $N=10$ which indicates that the error \blue{can remain small, even in higher dimensions.}

\section{Conclusion and future work}

In this paper, a novel shallow neural network is developed to approximate functions with jump discontinuities. The crucial idea is that a $d$-dimensional piecewise continuous function can be augmented as a continuous function defined in $(d+1)$-dimensional space. This function, based on the universal approximation theory, can be approximated by a shallow, feedforward, fully connected neural network. We thereby propose a simple neural network architecture consisting of an input layer of $(d+1)$-dimensional coordinate variables, a hidden layer with a moderate number of neurons, and an output layer of the function itself. We show that the present network is efficient and accurate for piecewise continuous function approximation, and can serve as a solution model for elliptic interface problems. Combined \blue{with the PINN-type loss functions}, the present network approximates the solution with high accuracy, and the results are comparable to the immersed interface method. Because of the mesh-free nature of the network, there is no difficulty in implementing the model for problems on irregular domains or in high dimensions.

The present work differs significantly from the work using DNNs in the literature as we only consider a completely shallow neural network (one hidden layer). To approximate functions or solutions that are highly oscillatory, one can increase the number of neurons in the hidden layer to achieve the desired accuracy. As we show in the numerical experiments in this paper, all the problems can be solved with high accuracy by the present network with a moderate number (less than a hundred) of neurons. Notice that a shallow network is much easier to train than a deep one.

In the present work, we only consider the stationary elliptic interface problems. As a forthcoming extension, we shall consider time-dependent problems, particularly the moving interface problems, which will be left as future work.

\section*{Acknowledgement}

W.-F. Hu, T.-S. Lin and M.-C. Lai acknowledge supports by the Ministry of Science and Technology,
Taiwan, under research grants 109-2115-M-008-014-MY2, 109-2115-M-009-006-MY2 and 110-2115-M-A49-011-MY3, respectively. The authors would like to thank Che-Chia Chang and Yi-Jun Shih for help with the development of the python LM code.



\end{document}